\documentclass
{amsart}
\usepackage{graphicx}
\newtheorem{theorem}{Theorem}[section]

\theoremstyle{definition}
\newtheorem{definition}[theorem]{Definition}
\newtheorem{example}[theorem]{Example}

\newcommand{\F}{\Bbb F}

\theoremstyle{remark}
\newtheorem{remark}[theorem]{Remark}

\numberwithin{equation}{section}

\begin{document}

\title[Reidemeister torsion and fibered knots]
{Reidemeister torsion, twisted Alexander polynomial 
and fibered knots} 

\author{
Hiroshi Goda, 
Teruaki Kitano 
and Takayuki Morifuji}

\thanks{The second and third author are supported 
in part by Grand-in-Aid for Scientific Research 
(No. 14740037 and No. 14740036 respectively), 
The Ministry of Education, Culture, 
Sports, Science and Technology, Japan.}

%\date{\today}

\begin{abstract}
As a generalization of a classical result on 
the Alexander polynomial for fibered knots, 
we show in this paper that the Reidemeister torsion 
associated to a certain representation detects 
fiberedness of knots in the three sphere.
\end{abstract}

\maketitle

\section{Introduction}

As is well-known, 
the Alexander polynomial of a fibered knot 
is monic (see \cite{Neu}, \cite{Rapa}, \cite{St}). 
That is, 
the coefficient of the highest degree term of 
the normalized Alexander polynomial is 
a unit $1\in\Bbb Z$. 
By the symmetry (or duality) 
of the Alexander polynomial, 
%we easily see that 
its lowest degree term is also one. 
This criterion is sufficient for 
alternating knots \cite{Murasugi} and 
prime knots up to 10 crossings \cite{kanenobu} 
for instance. 
However, in general, 
the converse is not true. 
In fact, 
there are infinitely many non-fibered knots 
having monic Alexander polynomials. 
If 
we remember here 
Milnor's result \cite{M1}, 
we have to remark that 
these claims on the Alexander polynomial 
can be restated by 
the Reidemeister torsion. 

%Recently, 
%a distance for links was 
%introduced by using 
%Heegaard splitting for sutured manifolds 
%(see \cite{G1},~\cite{G2}) and 
%multivalued Morse function 
%(or circle valued Morse function, 
%see \cite{N},~\cite{PRW}) 
%independently. 
%Roughly speaking, 
%it measures a distance between 
%a given link $L$ and the fibered links 
%by counting the number of 
%the $k$th link polynomial of $L$ 
%which is not monic (see \cite{PRW} for details). 
%It should be noted that 
%the Alexander polynomial is the 0th link polynomial \cite{Fox}
%and such an estimate is given by so-called 
%Morse-Novikov inequality. 
%From this viewpoint, 
%the above criterion using 
%monicness of the Alexander polynomial 
%seems to be very important. 

The purpose of this paper is 
to give a necessary condition that 
a knot in $S^3$ is fibered 
by virtue of the Reidemeister torsion 
associated to a certain linear representation. 
More precisely, 
we show that 
the Reidemeister torsion of a fibered knot 
defined for a certain tensor representation 
is expressed as a rational function of 
monic polynomials. 
This Reidemeister torsion is nothing but 
Wada's twisted Alexander polynomial 
(see \cite{Kitano} for details), 
so that 
our result can be regarded as 
a natural generalization of the property 
on the classical Alexander polynomial 
mentioned above. 
%Furthermore, 
%we expect that 
%the result obtained here 
%could contribute to get a sharp estimate 
%than the Morse-Novikov inequality. 

This paper is organized as follows. 
In the next section, 
we review the definition of Reidemeister torsion 
over a field $\F$. 
Further 
we describe how to compute it 
in the case of knot exteriors. 
The point of our method here is 
that the notion of monic makes sense 
for the Reidemeister torsion 
associated to a tensor representation of 
a unimodular representation over $\F$ 
and the abelianization homomorphism. 
In Section 3, 
we state and prove the main theorem of this paper. 
The final section 
is devoted to compute some examples. 

We should note here 
that there is a similar work by Cha \cite{C}. 
The notion of Alexander polynomials twisted by 
a representation and its applications have 
appeared in several papers 
(see \cite{JW},~\cite{KL},~\cite{Lin},~\cite{Wada}).

\section{Reidemeister torsion}

In this section, 
we review the definition of Reidemeister torsion 
over a field $\F$ (see \cite{J} and \cite{M2} for details). 

Let 
$V$ be an $n$-dimensional vector space over $\F$, 
and ${\bf b}=(b_1,\ldots,b_n)$ and ${\bf c}=(c_1,\ldots,c_n)$ 
two bases for $V$. 
If we put 
$\displaystyle{c_i=\sum_{j=1}^{n}a_{ij}b_j}$, 
we have a nonsingular matrix $A=(a_{ij})$ 
with coefficients in $\F$. 
Further 
let $[{\bf b}/{\bf c}]$ denote the determinant of $A$. 

Now let us consider an acyclic chain complex 
of finite dimensional vector spaces over $\F$: 
\[ 
C_*:
0 \longrightarrow C_m 
\overset{\partial_m}{\longrightarrow} C_{m-1}
\overset{\partial_{m-1}}{\longrightarrow}\cdots
\longrightarrow C_1
\overset{\partial_1}{\longrightarrow} C_0
\longrightarrow 0. 
\] 
We assume that 
a preferred basis ${\bf c}_q$ for $C_q(C_*)$ is given for any $q$. 
Choose any basis ${\bf b}_q$ of $B_q(C_*)$ and 
take its lift in $C_{q+1}(C_*)$, 
which we denote by $\widetilde{{\bf b}}_q$. 

Since the natural inclusion map 
$$
B_q(C_*)\to Z_q(C_*)
$$
is an isomorphism, 
the basis ${\bf b}_q$ can serve as a basis for $Z_q(C_*)$. 
Similarly 
the sequence 
\[
0 \longrightarrow Z_q(C_*)
  \longrightarrow C_q(C_*)
  \longrightarrow B_{q-1}(C_*)
  \longrightarrow 0
\]
is exact and the vectors 
$({\bf b}_q,\widetilde{{\bf b}}_{q-1})$ 
is a basis for $C_q(C_*)$. 
It is easily shown that 
$[{\bf b}_q,\widetilde{{\bf b}}_{q-1}/{\bf c}_{q}]$ 
is independent of 
the choices of $\widetilde{{\bf b}}_{q-1}$. 
Hence 
we may simply denote it by 
$[{\bf b}_q,{\bf b}_{q-1}/{\bf c}_{q}]$. 

\begin{definition}
The torsion of the chain complex $C_*$ is defined by 
the alternating product 
$$
\prod_{q=0}^{m}[{\bf b}_q,{\bf b}_{q-1}/{\bf c}_q]^{(-1)^{q+1}}
$$
and we denote it by $\tau(C_*)$.
\end{definition}

\begin{remark}
The torsion $\tau(C_*)$ depends only on the bases 
${\bf c}_0,\ldots,{\bf c}_m$.
\end{remark}

Now let us apply the above torsion 
to the following geometric situations. 
Let 
$X$ be a finite cell complex and 
$\widetilde{X}$ the universal covering of $X$ 
with the right action of $\pi_1X$ 
as deck transformations. 
Then 
the chain complex $C_*(\widetilde{X},\Bbb Z)$ 
has a structure of 
right free $\Bbb Z[\pi_1X]$-modules. 
Let 
$$
\rho:\pi_1X\to GL(n,\F)
$$
be a linear representation. 
We may regard $V$ as 
a $\pi_1X$-module by using $\rho$ 
and denote it by $V_\rho$. 
Define 
the chain complex $C_*(X,V_\rho)$ by 
$C_*(\widetilde{X},\Bbb Z)\otimes_{\Bbb Z[\pi_1X]}V_\rho$ 
and choose a preferred basis 
$$
\{
\sigma_1\otimes e_1,\sigma_1\otimes e_2,\ldots,
\sigma_1\otimes e_n,\ldots,\sigma_{k_q}\otimes e_1,\ldots,
\sigma_{k_q}\otimes e_n\}
$$
of $C_q(X,V_\rho)$, 
where $\{e_1,\ldots,e_n\}$ is a basis of $V$ and 
$\sigma_1,\ldots,\sigma_{k_q}$ are $q$-cells 
giving the preferred basis of $C_q(\widetilde{X},\Bbb Z)$.

Now 
we consider the following situation. 
That is, 
$C_*(X,V_\rho)$ is acyclic, 
in other words, all homology groups 
$H_*(X,V_\rho)$ vanish. 
In this case, 
we call $\rho$ an acyclic representation. 

\begin{definition}
Let 
$\rho:\pi_1X\to GL(n,\F)$ be an acyclic representation. 
Then 
Reidemeister torsion of $X$ with $V_\rho$-coefficients is 
defined by the torsion of the chain complex $C_*(X,V_\rho)$. 
We denote it by $\tau_\rho(X)$.
\end{definition}

\begin{remark}
It is known that 
$\tau_\rho(X)$ is well-defined as a PL-invariant, 
for an acyclic representation $\rho:\pi_1X\to GL(n,\F)$,
 up to a factor 
$\pm d$ where 
$d\in\mathrm{Im}(\mathrm{det}\circ\rho)\subset\F^*$.
As a reference, 
see \cite{M2} Section 8.
We can easily make a refinement of the above argument 
for our situation.
\end{remark}

Here 
let us consider a knot $K$ in $S^3$ and 
its exterior $E$. 
For 
the knot group $\pi K=\pi_1E$, 
we choose and fix a Wirtinger presentation 
$$
P(\pi K)=
\langle x_1,\ldots,x_u~|~r_1,\ldots,r_{u-1}\rangle.
$$
Then 
we can construct a 2-dimensional cell complex $X$ 
from $P(\pi K)$ such that 
$E$ collapses to $X$. 
The abelianization homomorphism 
$$
\alpha:\pi K\to H_1(E,\Bbb Z)\cong\Bbb Z
=\langle t\rangle
$$ 
is given by 
$$
\alpha(x_1)=\cdots=\alpha(x_u)=t.
$$
Furthermore, 
we always suppose that 
the image of a representation 
$\rho:\pi K\to GL(n,\F)$ 
is included in $SL(n,\F)$. 

These maps naturally induce the ring homomorphisms 
$\tilde{\rho}$ and $\tilde{\alpha}$ 
from $\Bbb Z[\pi K]$ to $M(n,\F)$ and $\Bbb Z[t^{\pm 1}]$ 
respectively, 
where 
$M(n,\F)$ denotes the matrix algebra of degree $n$ 
over $\F$. 
Then 
$\tilde{\rho}\otimes\tilde{\alpha}$ 
defines a ring homomorphism 
$$
\Bbb Z[\pi K]\to M\left(n,\F[t^{\pm 1}]\right).
$$
Let 
$F_u$ denote the free group with 
generators $x_1,\ldots,x_u$ and 
denote by 
$$
\Phi:\Bbb Z[F_u]\to M\left(n,\F[t^{\pm 1}]\right)
$$
the composite of the surjection 
$\Bbb Z[F_u]\to\Bbb Z[\pi K]$ 
induced by the presentation and the map 
$\Bbb Z[\pi K]\to M(n,\F[t^{\pm 1}])$ 
given by $\tilde{\rho}\otimes\tilde{\alpha}$. 

Let us consider the $(u-1)\times u$ matrix $M$ 
whose $(i,j)$ component is the $n\times n$ matrix 
$$
\Phi\left(\frac{\partial r_i}{\partial x_j}\right)
\in M\left(n,\F[t^{\pm 1}]\right),
$$
where 
$\displaystyle{\frac{\partial}{\partial x}}$ 
denotes the free differential calculus. 
This matrix $M$ is called the Alexander matrix of 
the presentation $P(\pi K)$ 
associated to the representation $\rho$. 

For 
$1\leq j\leq u$, 
let us denote by $M_j$ 
the $(u-1)\times(u-1)$ matrix obtained from $M$ 
by removing the $j$th column. 
We regard $M_j$ as 
a $(u-1)n\times (u-1)n$ matrix with coefficients in 
$\F[t^{\pm 1}]$. 

%On the other hand, 
Now let us recall that 
%we define 
the tensor representation 
$$
\rho\otimes\alpha:
\pi K\to GL(n,\F(t))
$$
is defined 
by 
$(\rho\otimes\alpha)(x)=\rho(x)\alpha(x)$ 
for $x\in \pi K$. 
Here 
$\F(t)$ denotes the rational function field over $\F$ 
and let $V$ be the $n$-dimensional vector space over $\F(t)$. 
Hereafter, 
%let us denote 
we denote 
the Reidemeister torsion $\tau_{\rho\otimes\alpha}(E)$ 
by $\tau_{\rho\otimes\alpha}K$. 

\begin{theorem}
All homology groups 
$H_*(E,V_{\rho\otimes\alpha})$ vanish 
$($namely, 
$\rho\otimes\alpha$ is an acyclic representation$)$ 
if and only if 
$\mathrm{det}~M_j\not=0$ for some $j$. 
In this case, 
we have 
$$
\tau_{\rho\otimes\alpha}K
=\frac{\mathrm{det}~M_j}{\mathrm{det}\Phi(x_j-1)},
$$
for any $j~(1\leq j\leq u)$. 
Moreover, 
$\tau_{\rho\otimes\alpha}$ is well-defined 
up to a factor 
$\pm t^{nk}~(k\in\Bbb Z)$ if $n$ is odd and 
up to only $t^{nk}$ if $n$ is even.
\end{theorem}

\begin{proof}
%It is easy to see that 
%${\mathrm{det}\Phi(x_j-1)}\neq 0$ 
%for any $i~(1\leq j\leq u)$ 
%because $\tilde{\alpha}(x_j-1)\neq 0$. 
%Then 
%it implies that 
%$H_0(E, V_{\rho\otimes \alpha})$ and 
%$H_1(E,V_{\rho\otimes \alpha})$ 
%are vanishing. 
%Since $E$ is a knot exterior, 
%its Euler number is zero. 
%Thus 
%the Euler number of $E$ with 
%$V_{\rho\otimes \alpha}$-coefficients 
%is also zero. 
%Since 
%$H_0$ and $H_1$ are vanishing, 
%$H_2(E,V_{\rho\otimes \alpha})$ is also vanishing. 
%
%Next 
%the above formula of $\tau_{\rho\otimes\alpha}$ 
%is obtained from easy calculations. 
%In fact, 
%we can reduce the computation of 
%$\tau_{\rho\otimes\alpha}K$ to the one of 
%$\tau_{\rho\otimes\alpha}(X)$ 
%by the simple homotopy invariance 
%of the Reidemeister torsion. 
%For example, see \cite{Kitano} Proposition 3.1 
%for more details. 
%In particular, 
%the ratio $\mathrm{det}~M_j/\mathrm{det}\Phi(x_j-1)$ 
%is independent of the choice of $j$ 
%(see \cite{Kitano} Lemma 1.2). 

The first two assertion are nothing but 
\cite{Kitano} Proposition 3.1. 
The independence on $j$ follows from 
\cite{Kitano} Lemma 1.2. 

%Finally, 
Next, 
if we consider well-definedness 
up to $\pm t^{nk}~(k\in\Bbb Z)$, 
we have only to recall Remark 2.4. 
The image of 
$$
\mathrm{det}\circ(\rho\otimes\alpha):
\pi K\to GL(n,\Bbb F(t))\to \Bbb F(t)^*
$$
is just $\{t^{nk}|k\in \Bbb Z\}$, 
because $\mathrm{Im}~\rho$ is 
included in $SL(n,\F)$. 
Therefore 
the claim immediately follows. 
Here 
if we take an even dimensional unimodular representation, 
$\tau_{\rho\otimes\alpha}$ 
is well-defined up to $t^{nk}$ 
(see also \cite{J} for details). 
\end{proof}

\begin{remark}
This theorem asserts that 
the twisted Alexander polynomial \cite{Wada} of a knot 
is the Reidemeister torsion of its knot exterior 
(see \cite{Kitano} for details). 
This is a generalization of Milnor's theorem in \cite{M1}. 
Recently 
this framework extended to more general situations by 
Kirk-Livingston in \cite{KL}. 
\end{remark}

\begin{remark}
Assume that 
$\rho$ is a homomorphism to $SL(n,R)$ 
over a unique factorization domain $R$ 
and the knot group $\pi K$ has a presentation 
which is strongly Tietze equivalent to a Wirtinger 
presentation of the knot. 
Then 
Wada shows in \cite{Wada} that 
the twisted Alexander polynomial of the knot 
associated to $\rho$ 
is well-defined up to 
a factor $\pm t^{nk}~(k\in\Bbb Z)$ 
if $n$ is odd and 
up to only $t^{nk}$ if $n$ is even.  
\end{remark}

\begin{remark}
If 
there is an element $\gamma$ of the commutator subgroup of $\pi K$ 
such that $1$ is not an eigenvalue of $\rho(\gamma)$, 
then $\tau_{\rho\otimes\alpha}$ becomes a ``polynomial" 
(see \cite{Wada}). 
Namely $\mathrm{det}~M_j$ is divided by $\mathrm{det}\Phi(x_j-1)$.
\end{remark}

\section{Main theorem}

In this section, 
we give a necessary condition that 
a knot $K$ in $S^3$ is fibered. 
A polynomial 
$a_mt^m+\cdots+ a_1t+a_0\in \F[t]$ 
is called monic 
if the coefficient $a_m$ is one. 
We then see from Theorem 2.5 that 
the notion of monic polynomial makes sense 
for the Reidemeister torsion. 

\begin{theorem}
For a fibered knot $K$ in $S^3$ 
and a unimodular representation $\rho:\pi K\to SL(2n,\F)$, 
the Reidemeister torsion $\tau_{\rho\otimes\alpha}K$ 
is expressed as a rational function of monic polynomials.
\end{theorem}

\begin{proof}
By using the fiber bundle structure, 
we can take the following presentation of $\pi K$:
$$
P(\pi K)=
\langle x_1,\ldots,x_{2g},h~|~
r_i=hx_ih^{-1}\varphi_*(x_i)^{-1},~1\leq i\leq 2g\rangle,
$$
where $x_1,\ldots,x_{2g}$ is a generating system 
of the fundamental group of the fiber surface of genus $g$, 
$h$ is a generator for $S^1$-direction 
corresponding to the meridian of $K$ 
and $\varphi_*$ denotes the automorphism of 
the surface group induced by the monodromy map $\varphi$. 
Here 
the abelianization homomorphism 
$\alpha:\pi K\to\Bbb Z=\langle t\rangle$ 
is given by 
$$
\alpha(x_1)=\cdots=\alpha(x_{2g})=1\quad \mathrm{and}\quad
\alpha(h)=t.
$$

This presentation of $\pi K$ allows us 
to define another 2-dimensional cell complex $Y$ 
instead of a cell complex $X$ constructed from 
a Wirtinger presentation of $\pi K$. 
Namely 
it has only one vertex, $2g+1$ edges, and $2g$ 2-cells 
attached by the relations of $P(\pi K)$.  
It is easy to see that there exists a homotopy equivalence 
$f:E\to Y$. 
From 
the result of Waldhausen \cite{Wal}, 
the Whitehead group $Wh(\pi K)$ is trivial 
for a knot group in general. 
Thus 
the Whitehead torsion of $f$ is also trivial element in $Wh(\pi K)$. 
It then implies that 
the homotopy equivalence map $f$ induces 
a simple homotopy equivalence from $E$ to $Y$. 
Since 
the Reidemeister torsion is a simple homotopy invariant, 
we can compute the Reidemeister torsion of $E$ 
as the one of $Y$ as follows. 
That is, 
we may use the previous presentation $P(\pi K)$ 
to compute $\tau_{\rho\otimes\alpha}K$ 
by means of Theorem 2.5. 

Let us consider the ``big" $2g\times 2g$ matrix $M$ 
whose $(i,j)$ component is the $n\times n$ matrix 
$$
\Phi\left(\frac{\partial r_i}{\partial x_j}\right)
\in M(n,\F[t^{\pm1}]).
$$
We then see that 
the diagonal component of $M$ is 
\begin{align*}
\Phi\left(\frac{\partial r_i}{\partial x_i}\right)
&=\Phi\left(h-\frac{\partial \varphi_*(x_i)}{\partial x_i}\right)\\
&=t\rho(h)
-\tilde{\rho}\left(\frac{\partial \varphi_*(x_i)}{\partial x_i}\right)
\end{align*}
and 
the coefficient of the highest degree term of 
$\mathrm{det}\Phi(\partial r_i/\partial x_i)$ is 
just $\mathrm{det}~\rho(h)=1$. 
Further 
other components $\Phi(\partial r_i/\partial x_j)~(i\not=j)$ 
do not contain $t$, 
so that 
the coefficient of the highest degree term of 
$\mathrm{det}~M$ is also one. 

On the other hand, 
the denominator is given by 
\begin{align*}
\mathrm{det}\Phi(h-1)
&= \mathrm{det}(t\rho(h)-I)\\
&=(\mathrm{det}~\rho(h))t^{2n}
-(\mathrm{tr}~\rho(h))t^{2n-1}+\cdots+1\\
&= t^{2n}+\cdots + 1,
\end{align*}
where 
$I$ denotes the identity matrix. 
Moreover 
$\rho$ is an even dimensional representation, 
so we see that $\tau_{\rho\otimes\alpha}K$ is 
well-defined up to a factor $t^{2nk}~(k\in\Bbb Z)$. 
This completes the proof. 
\end{proof}

\begin{remark}
If 
we can show directly that 
the presentation $P(\pi K)$ in Theorem 3.1 is 
strongly Tietze equivalent to 
a Wirtinger presentation of $K$, 
then 
Theorem 3.1 follows without using 
the result of Waldhausen (see Remark 2.7).
\end{remark}

\begin{remark}
If 
$\F$ is $\Bbb C$ or a finite field $\F_p$, 
then 
the Reidemeister torsion $\tau_{\rho\otimes\alpha}K$ 
for any knot $K$ and 
any representation 
$\rho:\pi K\to SL(2,\F)$ is symmetric. 
Namely, 
$\tau_{\rho\otimes\alpha}K$ is invariant 
under the transformation $t\mapsto t^{-1}$ 
up to a factor $t^k~(k\in \Bbb Z)$. 
Such a duality theorem appears 
originally in \cite{M1}. 
See also \cite{KL}, \cite{Kitano} and \cite{Mf} 
as for related works. 
\end{remark}

\section{Examples}

\begin{example}
Let 
$K$ be the figure eight knot $4_1$. 
This is one of the well-known genus one fibered knots 
in $S^3$. 
%Then 
The fundamental group of the exterior has 
a presentation 
$$
\pi K=
\langle
x,y~|~z x z^{-1}y^{-1}
\rangle,
$$
where $z=x^{-1}yxy^{-1}x^{-1}$. 
Let 
$\rho:\pi K\to SL(2,\Bbb C)$ be 
a noncommutative representation defined by 
$$
\rho(x)=\begin{pmatrix}1 & 1 \\ 0 &1\end{pmatrix}
\quad\mathrm{and}\quad
\rho(y)=\begin{pmatrix}1 & 0 \\ -\omega &1\end{pmatrix},
$$
where 
$\omega$ is a complex number satisfying 
$\omega^2+\omega+1=0$. 
As pointed out by Wada in \cite{Wada}, 
it is convenient to use relations 
instead of relators for 
the computation of the Alexander matrix. 
Thereby 
applying free differential calculus to 
the relation $r:zx=yz$, 
we obtain 
$$
\frac{\partial r}{\partial x}
=-x^{-1}+x^{-1}y+yx^{-1}-yx^{-1}y+yx^{-1}yxy^{-1}x^{-1}.
$$
Thus 
we have the matrix 
$$
M_2
=\left(\Phi\left(\frac{\partial r}{\partial x}\right)\right)
=
\begin{pmatrix}
-(\omega+1)t+\omega+2-t^{-1} & t+\omega-2+t^{-1} \\
(\omega-1)t-\omega+1 & -(\omega+1)t+3-t^{-1}\end{pmatrix}.
$$
Then 
the numerator of $\tau_{\rho\otimes\alpha}$ is given by 
\begin{align*}
\mathrm{det}~M_2
&= t^{-2}(t^4-6t^3+\omega^4t^2+\omega^2t^2+11t^2-6t+1)\\
&= t^{-2}(t-1)^2(t^2-4t+1).
\end{align*}
On the other hand, 
the denominator of $\tau_{\rho\otimes\alpha}$ 
is given by 
\begin{align*}
\mathrm{det}\Phi(y-1)
&=\mathrm{det}(t\rho(y)-I)\\
&=\mathrm{det}\begin{pmatrix}t-1&0\\-\omega t&t-1\end{pmatrix}\\
&=(t-1)^2.
\end{align*}
Therefore 
the Reidemeister torsion of the figure eight knot $K$ is 
$$
\tau_{\rho\otimes\alpha}K
=t^2-4t+1
$$
and this is in fact a monic polynomial.
\end{example}

\begin{example}
Let $KT$ be the Kinoshita-Terasaka knot \cite{KT}. 
It is well-known that 
$KT$ is one of the classical examples of knots 
with the trivial Alexander polynomial. 
The knot group $\pi KT$ has 
a presentation with four generators 
$x_1,\ldots,x_4$ and three relations \cite{Wada}: 
\begin{align*}
r_1&: x_1x_2x_1^{-1} = x_4x_2x_4x_2^{-1}x_4^{-1},\\
r_2&: x_4x_2x_4^{-1} 
= x_2^{-1}x_3x_1x_3^{-1}x_2x_1x_2^{-1}x_3x_1^{-1}x_3^{-1}x_2,\\
r_3&: x_1x_3x_1^{-1} = x_4x_3x_4x_3^{-1}x_4^{-1}.
\end{align*}
Applying free differential calculus, 
we have 
$$
\frac{\partial r_1}{\partial x_1}=1-x_1x_2x_1^{-1},\qquad
\frac{\partial r_1}{\partial x_2}=x_1-x_4+x_4x_2x_4x_2^{-1},\qquad
\frac{\partial r_1}{\partial x_3}=0,
$$
\begin{align*}
\frac{\partial r_2}{\partial x_1}=
& - x_2^{-1}x_3 - x_2^{-1}x_3x_1x_3^{-1}x_2 
+ x_2^{-1}x_3x_1x_3^{-1}x_2x_1x_2^{-1}x_3x_1^{-1},\\
\frac{\partial r_2}{\partial x_2}=
& ~x_4+x_2^{-1}
-x_2^{-1}x_3x_1x_3^{-1}
+ x_2^{-1}x_3x_1x_3^{-1}x_2x_1x_2^{-1}\\
& - x_2^{-1}x_3x_1x_3^{-1}x_2x_1x_2^{-1}x_3x_1^{-1}x_3^{-1},\\
\frac{\partial r_2}{\partial x_3}= 
& - x_2^{-1}
+ x_2^{-1}x_3x_1x_3^{-1}
- x_2^{-1}x_3x_1x_3^{-1}x_2x_1x_2^{-1}\\
& + x_2^{-1}x_3x_1x_3^{-1}x_2x_1x_2^{-1}x_3x_1^{-1}x_3^{-1},
\end{align*}
$$
\frac{\partial r_3}{\partial x_1}=1-x_1x_3x_1^{-1},\qquad
\frac{\partial r_3}{\partial x_2}=0,\qquad
\frac{\partial r_3}{\partial x_3}=x_1-x_4+x_4x_3x_4x_3^{-1}.
$$

Let 
$\rho:\pi KT\to SL(2,\F_5)$ be a 
noncommutative representation over 
the finite field $\F_5$ defined as follows:
$$
\rho(x_1)=\begin{pmatrix}0&1\\4&1\end{pmatrix},\quad
\rho(x_2)=\begin{pmatrix}0&4\\1&1\end{pmatrix},\quad
\rho(x_3)=\begin{pmatrix}0&1\\4&1\end{pmatrix}\quad\mathrm{and}\quad
\rho(x_4)=\begin{pmatrix}4&4\\3&2\end{pmatrix}.
$$
Then 
we obtain 
$$
M_4=
\begin{pmatrix}
3t+1&t&t&t^2+2t&0&0\\
2t&t+1&4t^2+t&4t&0&0\\
1&4t+3&t+2+t^{-1}&t+t^{-1}&3t+3+4t^{-1}&3t+4t^{-1}\\
4t&2t+1&t+1+4t^{-1}&3t+4&2t+4+t^{-1}&4t+1\\
1&4t&0&0&3t^2+t&2t^2+2t\\
t&4t+1&0&0&4t^2+t&3t^2+4t
\end{pmatrix}.
$$
Therefore 
the Reidemeister torsion of $KT$ is given by
\begin{align*}
\tau_{\rho\otimes\alpha}KT
&=\frac{\mathrm{det}~M_4}{\mathrm{det}\Phi(x_4-1)}\\
&=\frac{t^2(4t^8+t^7+t^6+4t^5+3t^4+4t^3+t^2+t+4)}{t^2+4t+1}\\
&=4t^6 + 2t^4 + t^3 + 2t^2 + 4.
\end{align*}
This is well-defined up to a factor $t^{2k}\ (k\in \Bbb Z)$, 
so that 
we 
%can 
may conclude 
the Kinoshita-Terasaka knot $KT$ 
is not fibered.
\end{example}

\begin{example}
Let $K$ be the knot illustrated in Figure 1. 
The normalized Alexander polynomial of $K$ 
is equal to the monic polynomial 
$t^4-t^3+t^2-t+1$. 
The knot group $\pi K$ has a presentation 
with seven generators $x_1, \ldots, x_7$ and 
six relations: 
\begin{align*}
r_1&: x_2x_1=x_3x_2x_1x_2x_1^{-1}x_2^{-1},\\
r_2&: x_6x_5x_6^{-1}=x_4x_3x_1^{-1}x_3x_1^{-1}x_3x_1x_3^{-1}x_1x_3^{-1}x_1x_3^{-1}x_4^{-1},\\
r_3&: x_6x_7x_6^{-1}=x_4x_3x_1^{-1}x_3x_1^{-1}x_3x_1x_3^{-1}x_1x_3^{-1}x_4^{-1} ,\\
r_4&: x_5x_6x_5^{-1}=x_7x_2x_7^{-1} ,\\
r_5&: x_2x_6x_2^{-1}=x_3x_2x_1x_2x_1^{-1}x_2^{-1}x_3^{-1}x_7x_3x_2x_1x_2^{-1}x_1^{-1}x_2^{-1}x_3^{-1},\\
r_6&: x_5x_4x_5^{-1}x_7=x_7x_3x_2x_1x_2x_1^{-1}x_2^{-1}x_3^{-1}.
\end{align*}

Let 
$\rho:\pi K\to SL(2, \F_5)$ be a 
noncommutative representation over 
$\F_5$ defined as follows:
$$
\rho(x_1)=\begin{pmatrix}1&1\\0&1\end{pmatrix},\quad
\rho(x_2)=\begin{pmatrix}1&0\\4&1\end{pmatrix},\quad
\rho(x_3)=\begin{pmatrix}1&0\\4&1\end{pmatrix},\quad
\rho(x_4)=\begin{pmatrix}2&1\\4&0\end{pmatrix},\quad$$
$$\rho(x_5)=\begin{pmatrix}2&4\\1&0\end{pmatrix},\quad
\rho(x_6)=\begin{pmatrix}3&1\\1&4\end{pmatrix}\quad
\mathrm{and}\quad
\rho(x_7)=\begin{pmatrix}1&4\\0&1\end{pmatrix}.
$$

By the same method as in previous examples, 
we have the following Reidemeister torsion of $K$:
\begin{align*}
\tau_{\rho\otimes\alpha}K
&=\frac{\mathrm{det}~M_7}{\mathrm{det}\Phi(x_7-1)}\\
&=\frac{t^{12}(3t^4+4t^3+t^2+4t+3)}{t^2+3t+1}\\
&=3t^2+3.
\end{align*}
Hence this knot $K$ is not fibered.
\end{example}

\begin{figure}
\centering
\includegraphics[width=7.2cm,height=3.3cm]{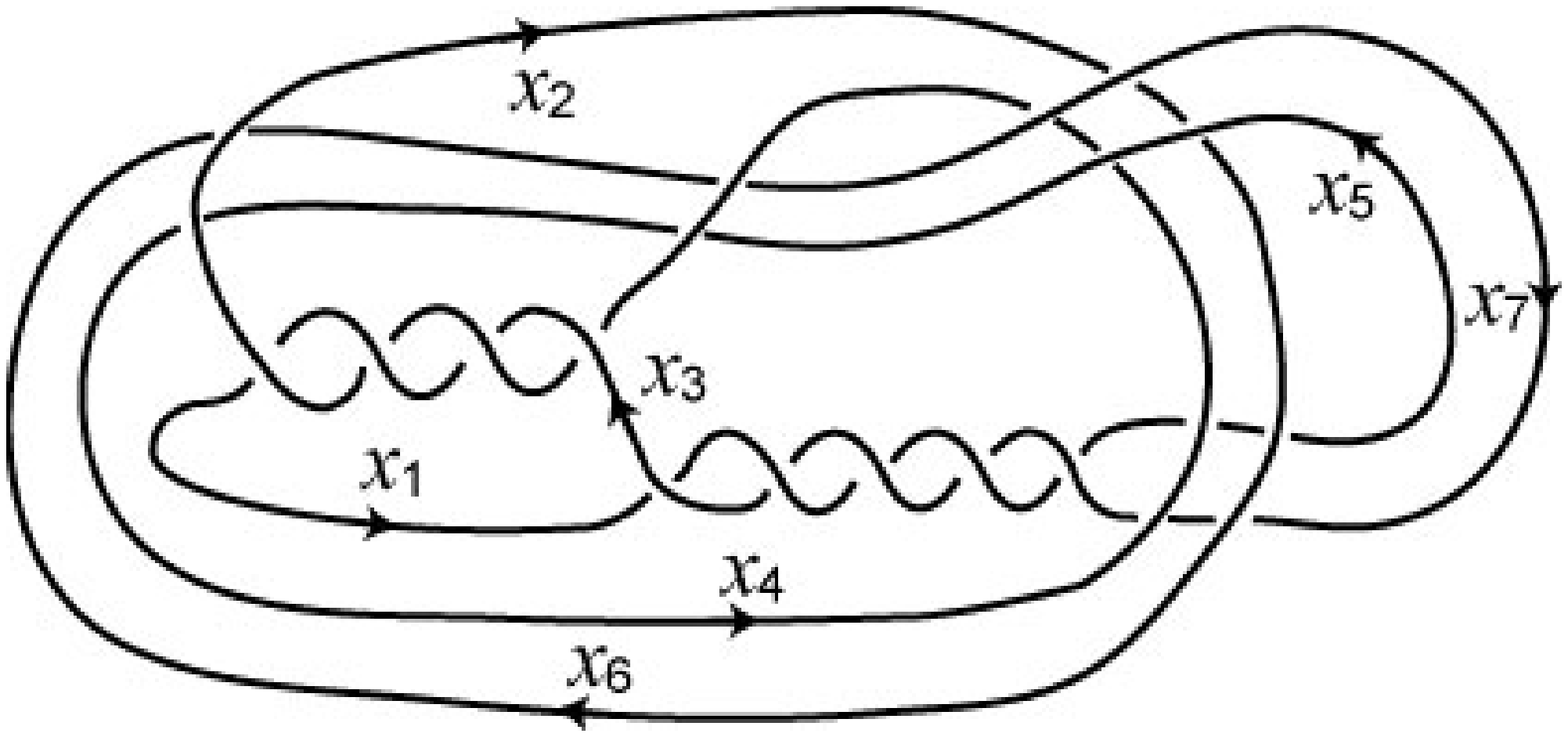}
\caption{}
\end{figure}

%\begin{remark}
We use Kodama's program ``KNOT'' and Wada's one 
to compute these examples. 
The former is to obtain $N$-data (see \cite{Wada0}) 
from a knot projection, 
which are necessary 
to input into Wada's program. 
The latter one is to 
compute unimodular representations over finite fields 
of knot groups from $N$-data. 
Here 
we should remark that 
Kodama's one works on Linux while 
Wada's one works on Macintosh. 
%If both of them are amalgamated, 
%then 
%it could be possible to detect fiberedness 
%of a knot from the projection in a few seconds
%almost completely. 
%\end{remark}

\vspace{2mm}

\noindent
{\it Acknowledgements}. 
The authors would like to thank 
Dr. Mitsuhiko Takasawa and Dr. Masaaki Suzuki 
for their useful advices on the computer. 
The authors would like to express their 
thanks to Professor Andrei Pajitnov for his 
careful reading the earlier version of 
this paper. 

\

\bibliographystyle{amsplain}

\begin{thebibliography}{30}

\bibitem{C} J.C. Cha,
\textit{Fibered knots and twisted Alexander invariants}, 
preprint. 

%\bibitem{Fox}
%R. Crowell and R. Fox,  
%Introduction to knot theory,
%Ginn and Co., Boston, Mass. (1963). 


%\bibitem{G1}
%H. Goda,
%\textit{Heegaard splitting for sutured manifolds and Murasugi sum},
%Osaka J. Math. 29 (1992), 21-40.

%\bibitem{G2}
%H. Goda,
%\textit{On handle number of Seifert surfaces in $S^{3}$}, 
%Osaka J. Math. 30 (1993), 63-80.

\bibitem{JW}
B. Jiang and S. Wang, 
\textit{Twisted topological invariants associated with representations}, 
in Topics in Knot Theory (1993), 211--227. 

\bibitem{J}D. Johnson, 
\textit{A geometric form of Casson's invariant and its 
connection to Reidemeister torsion}, 
unpublished lecture notes. 

\bibitem{kanenobu}
T. Kanenobu, 
\textit{The augmentation subgroup of a pretzel link}, 
Math. Sem. Notes Kobe Univ. 7 (1979), 363--384.

\bibitem{KT} 
S. Kinoshita and H. Terasaka, 
\textit{On unions of knots}, 
Osaka Math. J. 9 (1957), 131-153. 

\bibitem{KL}P. Kirk and C. Livingston, 
\textit{Twisted Alexander invariants, Reidemeister torsion, 
and Casson-Gordon invariants}, 
Topology 38 (1999), 635--661.

\bibitem{Kitano} T. Kitano,
\textit{Twisted Alexander polynomial and Reidemeister torsion},
Pacific J. Math. 174 (1996), 431--442

\bibitem{Lin}
X.S. Lin, 
\textit{Representations of knot groups and twisted Alexander polynomials}, 
Acta Math. Sin. (Engl. Ser.) 17 (2001), 361--380.

\bibitem{M1}J. Milnor, 
\textit{A duality theorem for Reidemeister torsion}, 
Ann. of Math. 76 (1962), 137--147.

\bibitem{M2} J. Milnor, 
\textit{Whitehead torsion},  
Bull. Amer. Math. Soc. 72 (1966), 358--426.

\bibitem{Mf}T. Morifuji, 
\textit{Twisted Alexander polynomial for the braid group}, 
Bull. Austral. Math. Soc. 64 (2001), 1--13.


\bibitem{Murasugi}
K. Murasugi,
\textit{On a certain subgroup of the group of an alternating link},
Amer. J. Math. 85 (1963), 544--550.



\bibitem{Neu}
L. Neuwirth,
Knot Groups, 
Annals of Mathematics Studies, No. 56 
Princeton University Press, Princeton, N.J.(1965). 

%\bibitem{N}
%S. P. Novikov, 
%\textit{Multivalued functions and functionals, 
%An analogue of the Morse theory}, 
%Soviet Math.Dokl.24 (1981), 222-226. 



%\bibitem{PRW}
%A. Pajitnov, L. Rudolph, and C. Weber, 
%\textit{Morse-Novikov number for knots and links}, 
%St. Petersburg Math. J. 13 (2002), 417-426.



\bibitem{Rapa}
E. Rapaport,
\textit{On the commutator subgroup of a knot group}, 
Ann. of Math. (2) 71 (1960), 157--162


\bibitem{S} L.C. Siebenmann, 
\textit{A total Whitehead torsion obstruction to fibering over the circle}, 
Comment. Math. Helv. 45 (1970),  1--48. 

\bibitem{St}
J. Stallings, 
\textit{On fibering certain $3$-manifolds}, 
1962 Topology of 3-manifolds and related topics 
(Proc. The Univ. of Georgia Institute, 1961), 95--100.


\bibitem{Wada0}
M. Wada, 
\textit{Coding link diagrams},  
J. Knot Theory Ramifications 2 (1993), 233--237.

\bibitem{Wada} M. Wada,
\textit{Twisted Alexander polynomial for finitely 
presentable groups},
Topology 33 (1994), 241--256.


\bibitem{Wal} F. Waldhausen,
\textit{Algebraic $K$-theory of generalized free products. I, II},
Ann. of Math. (2) 108 (1978),  135--204.

\end{thebibliography}

\

\begin{small}

\begin{flushleft}
Hiroshi Goda: goda@cc.tuat.ac.jp \\
Department of Mathematics \\ 
Tokyo University of Agriculture and Technology \\ 
2-24-16 Naka-cho, Koganei \\ 
Tokyo 184-8588, Japan \\ 

\

Teruaki Kitano: kitano@is.titech.ac.jp \\
Department of Mathematical and Computing Sciences \\ 
Tokyo Institute of Technology \\ 
2-12-1 Oh-okayama, Meguro-ku \\ 
Tokyo 152-8552, Japan \\ 

\

Takayuki Morifuji: morifuji@cc.tuat.ac.jp \\
Department of Mathematics \\ 
Tokyo University of Agriculture and Technology \\ 
2-24-16 Naka-cho, Koganei \\ 
Tokyo 184-8588, Japan \\ 
\end{flushleft}
\end{small}

\end{document}